# KARL PEARSON'S META-ANALYSIS REVISITED

By Art B. Owen[1]

*Stanford University*

This paper revisits a meta-analysis method proposed by Pearson [*Biometrika* **26** (1934) 425–442] and first used by David [*Biometrika* **26** (1934) 1–11]. It was thought to be inadmissible for over fifty years, dating back to a paper of Birnbaum [*J. Amer. Statist. Assoc.* **49** (1954) 559–574]. It turns out that the method Birnbaum analyzed is not the one that Pearson proposed. We show that Pearson's proposal is admissible. Because it is admissible, it has better power than the standard test of Fisher [*Statistical Methods for Research Workers* (1932) Oliver and Boyd] at some alternatives, and worse power at others. Pearson's method has the advantage when all or most of the nonzero parameters share the same sign. Pearson's test has proved useful in a genomic setting, screening for age-related genes. This paper also presents an FFT-based method for getting hard upper and lower bounds on the CDF of a sum of nonnegative random variables.

**1. Introduction.** Methods for combining $p$-values have long been studied. Recent research in genomics [Zahn et al. (2007)] and functional magnetic resonance imaging (fMRI) [see references in Benjamini and Heller (2007)] has sparked renewed interest. One gets a large matrix of $p$-values and considers meta-analysis along the rows, controlling for multiplicity issues in the resulting column of combined $p$-values.

This paper revisits an old issue in meta-analysis. A genomic application lead to the reinvention of a meta-analysis technique of Pearson (1934). That method has long been out of favor because the paper by Birnbaum (1954) appears to have proved that it is inadmissible (in an exponential family context). This paper shows that the method is in fact admissible in that exponential family context. Being admissible, it is more powerful than the widely used combination method of Fisher (1932) at some points of the

Received August 2008; revised January 2009.
[1]Supported by Grant DMS-06-04939 from the U.S. National Science Foundation.
*[AMS 2000 subject classifications](http://www.ams.org/msc/).* 62F03, 62C15, 65T99, 52A01.
*Key words and phrases.* Admissibility, fast Fourier transform, hypothesis testing, microarrays.

alternative hypothesis. Such points turn out to be especially important in the motivating biological problem.

The motivating work, reported in the AGEMAP study of Zahn et al. (2007), was to screen 8932 genes, searching for those having expression levels that are correlated with age in the mouse. Each gene was tested in $m = 16$ tissues, yielding 16 regression coefficients. There was an $8932 \times 16$ matrix of $p$-values. The data were noisy and it was plausible a priori that most if not all genes could have little or no correlation with age. In addition to the recently well-studied issue of controlling for multiple testing over genes, there was a problem of pooling information from different tissues for any single gene. For issues of multiple testing see Dudoit and van der Laan (2008). This article focuses on the problem of basing a single test on $m > 1$ test statistics.

A gene that is not age-related has a slope parameter of zero in all $m$ tissues. For a gene that is age-related, we expect several nonzero slopes. Because gene expression is tissue-dependent, it is also possible that a gene correlated with age in numerous tissues might fail to correlate in some others. Therefore, assuming a common nonzero slope is unreasonable. It is even possible that a gene's expression could increase with age in some tissues while decreasing in one or more others. But we do expect that the nonzero slopes should be predominantly of the same sign. The prior understanding of the biology is not detailed enough to let us specify in advance how many tissues might be nonage-related or even discordant for an otherwise age-related tissue.

Pearson's method is well suited to this problem. The better-known method of Fisher combines $p$-values by taking their product. It can be based on one-tailed or two-tailed test statistics. When based on one-tailed statistics, Fisher's test works best if we know the signs of the alternatives. When based on two-tailed statistics, Fisher's test does not favor alternatives that share a common sign. The test of Pearson (1934) may be simply described. It runs a Fisher-style test for common left-sided alternatives as well as one for common right-sided alternatives, and it takes whichever of those two is most extreme. It builds in a strong preference for common directionality while not requiring us to know the common direction, and it remains powerful when a small number of tests differ in sign from the dominant one. These properties fit the needs in Zahn et al. (2007).

An outline of this paper is as follows. Section 2 defines the test statistics and hypotheses that we work with. Section 3 reviews some basic concepts in meta-analysis and compares Pearson's test graphically to the better-known competitors. Section 4 shows that Pearson's method is admissible in the exponential family context. It also shows that a simple factor of two Bonferroni correction is very accurate for tail probabilities based on Pearson's test statistic. Section 5 reviews the history surrounding the misinterpretation of Pearson's proposal. Part of the problem stems from mixing up $p$-values from one and two-tailed tests. Section 6 compares the power of Pearson's test with some others, including tests that are based on the original test statistics, not just the $p$-values. Section 7 considers some of the recent literature on combining $p$-values. A discussion is in Section 8. The Appendix presents the numerical methods used to make the power computations. These are based on a kind of interval arithmetic for sub-distributions using the FFT.

**2. Notation.** This section sets out the notation for the paper. First are the parameters, hypotheses and test statistics for the univariate settings. Then they are combined to define multivariate rejection regions and test statistics.

2.1. *Parameters and hypotheses.* We consider a setting where there are $m$ parameters $\beta_1, \ldots, \beta_m$ and $m$ corresponding estimates $\hat{\beta}_1, \ldots, \hat{\beta}_m$. These estimates are random variables whose observed values are denoted by $\hat{\beta}_1^{\mathrm{obs}}, \ldots, \hat{\beta}_m^{\mathrm{obs}}$. In the motivating problem, these were regression slopes. We assume that the $m$ statistics $\hat{\beta}_j$ are independent. Dependent test statistics are considered briefly in Section 6.3. We also assume that $\hat{\beta}_j$ have continuous distributions.

For $j = 1, \ldots, m$, we can consider the hypotheses

$$H_{0,j} : \beta_j = 0$$
$$H_{L,j} : \beta_j < 0$$
$$H_{R,j} : \beta_j > 0$$

and

$$H_{U,j} : \beta_j \neq 0,$$

based on the sign of $\beta_j$. These are the null hypotheses, left- and right-sided alternatives and an undirected alternative, respectively.

Using $\hat{\beta}_j^{\mathrm{obs}}$ as test statistics, we may define

$$\widetilde{p}_j = \Pr(\hat{\beta}_j \leq \hat{\beta}_j^{\mathrm{obs}} | \beta_j = 0)$$

and

$$p_j = \Pr(|\hat{\beta}_j| \geq |\hat{\beta}_j^{\mathrm{obs}}| | \beta_j = 0).$$

The $p$-values for alternatives $H_{L,j}$, $H_{R,j}$ and $H_{U,j}$, respectively, are $\widetilde{p}_j$, $1 - \widetilde{p}_j$ and $p_j = 2\min(\widetilde{p}_j, 1 - \widetilde{p}_j)$. A $p$-value for a one-tailed or two-tailed test is called a one-tailed or two-tailed $p$-value below.

For the entire parameter vector $\beta = (\beta_1, \ldots, \beta_m)$, it is straightforward to define the simple null hypotheses $H_0$ for which $\beta_1 = \beta_2 = \cdots = \beta_m = 0$. We do not consider composite null hypotheses.

For $m > 1$, the alternatives to $H_0$ are more complicated than for $m = 1$. There are $3^m$ possible values for the vector of signs of $\beta_j$ values and many possible subsets of these could be used to define the alternative. For example, one could take

$$H_L = (-\infty, 0]^m - \{0\} \quad \text{or} \quad H_R = [0, \infty)^m - \{0\},$$

or their union.

But any of these choices leaves out possibilities of interest. Therefore, we take the alternative $H_A$ to be that $\beta_j \neq 0$ for at least one $j \in \{1, \ldots, m\}$. That is, $H_A : \beta \in \mathbb{R}^m - \{0\}$.

While all of $\mathbb{R}^m - \{0\}$ is of interest, the parts with concordant signs are of greater interest than those with discordant signs. For example, with $\Delta > 0$, we want greater power against alternatives $\pm(\Delta, \Delta, \ldots, \Delta)$ than against other alternatives of the form $(\pm\Delta, \pm\Delta, \ldots, \pm\Delta)$. In a screening problem, the former are more convincing, while the latter cause one to worry about noise and systematic experimental biases. The situation is analogous to the choice of the Tukey's versus Sheffé's statistic in multiple comparisons: both have the alternative of unequal means, but their power versus specific alternatives of interest could lead us to prefer one to the other in a given application.

It is worthwhile to represent the vector $\beta$ as the product $\tau\theta$, where $\theta \in \mathbb{R}^m$ is a unit vector, and $\tau \geq 0$. We may then consider the power of various tests of $H_0$ as $\tau$ increases.

2.2. *Rejection regions.* The decision to accept or reject $H_0$ will be based on $\widetilde{p}_1, \ldots, \widetilde{p}_m$. As usual, acceptance really means failure to reject and is not interpreted as establishing $H_0$. The rejection region is $\widetilde{R} = \{(\widetilde{p}_1, \ldots, \widetilde{p}_m) | H_0 \text{ rejected}\} \subset [0, 1]^m$. For some of the methods we consider, the rejection region can be expressed in terms of the two-tailed $p$-values $p_j$. Then we write $R = \{(p_1, \ldots, p_m) | H_0 \text{ rejected}\} \subset [0, 1]^m$.

While formulating the region in terms of $p$-values seems unnatural when the raw data are available, it does not change the problem much. For example, if $\hat{\beta}_j \sim \mathcal{N}(\beta_j, \sigma_j^2)$ with known $\sigma_j^2$, then any region defined through $\widetilde{p}_j$-values can be translated into one for $\beta_j$-values. In that case, we write $\widetilde{R}' = \{(\hat{\beta}_1, \ldots, \hat{\beta}_m) | H_0 \text{ rejected}\}$. It is more realistic for the $p$-values to come from a $t$ distribution based on estimates of $\sigma_j^2$. For large degrees of freedom, the normal approximation to the $t$ distributed problem is a reasonable one, and it is simpler to study. Discussion of $t$ distributed test statistics is taken up briefly in Section 6.1.

2.3. *Test statistics.* Under the null hypothesis $\widetilde{p}_j$, $1 - \widetilde{p}_j$ and $p_j$ all have the $U(0, 1)$ distribution. It follows that

$$(2.1) \qquad Q^L \equiv -2\log\left(\prod_{j=1}^{m} \widetilde{p}_j\right),$$

$$(2.2) \quad Q^R \equiv -2\log\left(\prod_{j=1}^{m}(1-\widetilde{p}_j)\right) \quad \text{and}$$

$$(2.3) \quad Q^U \equiv -2\log\left(\prod_{j=1}^{m}p_j\right)$$

all have the $\chi^2_{(2m)}$ distribution under $H_0$. An $\alpha$ level test based on any of these quantities rejects $H_0$ when that $Q$ is greater than or equal to $\chi^{2,1-\alpha}_{(2m)}$. Their chi-square distribution is due to Fisher (1932).

When $m=1$, these three tests reduce to the usual one and two-sided tests. When $m>1$ they are reasonable generalizations of one and two-sided tests.

The test of Pearson (1934) is based on

$$(2.4) \quad Q^C \equiv \max(Q^L, Q^R).$$

If $m=1$, then $Q^C = Q^U$, but for $m>1$ they differ. The superscript $C$ is mnemonic for concordant.

The null distribution of $Q^C$ is not $\chi^2_{(2m)}$. However, a Bonferroni correction is quite accurate:

$$(2.5) \quad \alpha - \frac{\alpha^2}{4} \leq \Pr(Q^C \geq \chi^{2,1-\alpha/2}_{(2m)}) \leq \alpha.$$

Equation (2.5) follows from Corollary 1 in Section 4. For instance, when the nominal level is $\alpha = 0.01$, the attained level is between 0.01 and 0.009975. Equation (2.5) shows that $\min(1, 2\Pr(\chi^2_{(2m)} \geq Q^C))$ is a conservative $p$-value. Equation (2.5) shows that the accuracy of this Bonferroni inequality improves for small $\alpha$ which is where we need it most.

The statistic $Q^L$ is the natural one to use when the alternative is known to be in $H_L$. But it still has power tending to 1, as $\tau = \|\beta\|$ tends to infinity, so long as $\theta = \beta/\tau$ is not in $H_R$. Naturally, $Q^R$ has a similar problem in $H_L$ while being well-suited for $H_R$. Neither $Q^C$ nor $Q^U$ have such problems. If we have no idea which orthant $\beta$ might be in, then $Q^U$ is a natural choice, while if we suspect that the signs of large (in absolute value) nonzero $\beta_j$ are mostly the same, then $Q^C$ has an advantage over $Q^U$.

2.4. *Stouffer et al.'s meta-analysis.* An alternative to Fisher's method is that of Stouffer et al. (1949), which is based on turning the $p$-values into $Z$-scores. Let $\varphi(x) = \exp(-x^2/2)/\sqrt{2\pi}$ denote the $\mathcal{N}(0,1)$ density, and then let $\Phi(x) = \int_{-\infty}^{x} \varphi(z)\,dz$.

We can define tests of $H_0$ based on $Z$-scores via

$$(2.6) \quad S^L = \frac{1}{\sqrt{m}} \sum_{j=1}^{m} \Phi^{-1}(1-\widetilde{p}_j),$$

$$(2.7) \qquad S^R = \frac{1}{\sqrt{m}} \sum_{j=1}^{m} \Phi^{-1}(\widetilde{p}_j),$$

$$(2.8) \qquad S^U = \frac{1}{\sqrt{m}} \sum_{j=1}^{m} |\Phi^{-1}(\widetilde{p}_j)| \quad \text{and}$$

$$(2.9) \qquad S^C = \max(S^L, S^R),$$

which are directly analogous to the tests of Section 2.3. For independent tests considered here, $S^L$ and $S^R$ have the $\mathcal{N}(0,1)$ distribution under $H_0$, while $S^C$ has a half-normal distribution and $S^U$ does not have a simple distribution. Note that $S^L = -S^R$ and that $S^C = |S^L| = |S^R|$.

**3. Meta-analysis and a graphical comparison of the tests.** This section reviews some basics of meta-analysis for further use. Then it presents a graphical comparison of Pearson's test with the usual tests, to show how it favors alternatives with concordant signs. For background on meta-analysis, see Hedges and Olkin (1985).

It has been known since Birnbaum (1954) that there is no single best combination of $m$ independent $p$-values. A very natural requirement for a combination test is Birnbaum's.

CONDITION 1. If $H_0$ is rejected for any given $(p_1, \ldots, p_m)$, then it will also be rejected for all $(p_1^*, \ldots, p_m^*)$ such that $p_j^* \leq p_j$ for $j = 1, \ldots, m$.

Birnbaum proved that every combination procedure which satisfies Condition 1 is in fact optimal, for some monotone alternative distribution. Optimality means maximizing the probability of rejecting $H_0$, subject to a constraint on the volume of the region $R$ of vectors $(p_1, \ldots, p_m)$, for which $H_0$ is rejected. Birnbaum allows simple alternatives that have independent $p_j$ with decreasing densities $g_j(p_j)$ on $0 \leq p_j \leq 1$. He also allows Bayes mixtures of such simple alternatives. Birnbaum shows that Condition 1 is necessary and sufficient for admissibility of the combination test, again in the context of decreasing densities. Here is his definition of admissibility:

DEFINITION 1 [Birnbaum (1954), page 564]. A test is admissible if there is no other test with the same significance level, which, without ever being less sensitive to possible alternative hypotheses, is more sensitive to at least one alternative.

The top row of Figure 1 illustrates 4 rejection regions $R$ satisfying Condition 1, arising from the Fisher test, the Stouffer test, a test based on $\min(p_1, \ldots, p_m)$ and a test based on $\max(p_1, \ldots, p_m)$, for the case $m = 2$.

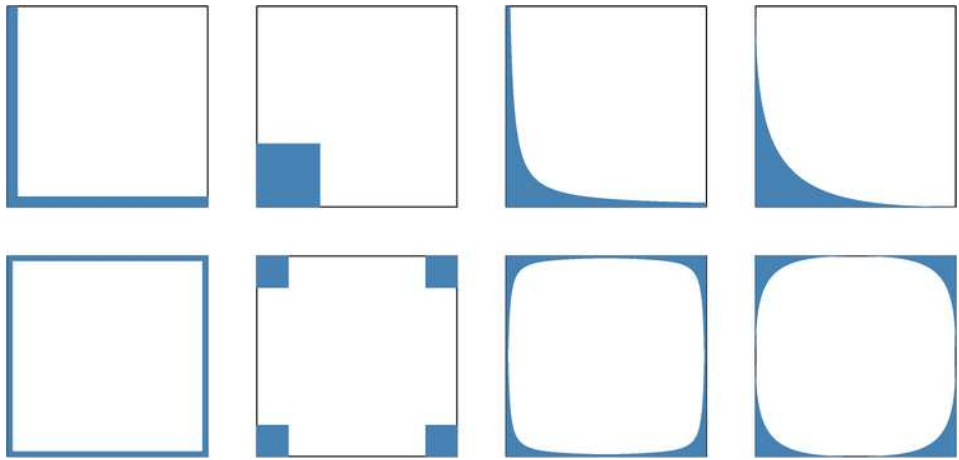

FIG. 1. *This figure shows rejection regions for a pair of tests. The top four images have coordinates $(p_1, p_2)$ where $p_j$ near zero is evidence against $H_{0j}$ in a two-tailed test. The columns, from left to right, are based on $\min(p_1, p_2)$, $\max(p_1, p_2)$, Fisher's combination and Stouffer's combination, as described in the text. Each region has area $1/10$. The bottom row shows the same rejection regions in coordinates $(\widetilde{p}_1, \widetilde{p}_2)$, where $\widetilde{p}_j$ near 0 is evidence that $\beta_j < 0$, and $\widetilde{p}_j$ near 1 is evidence that $\beta_j > 0$.*

Using the minimum is due to Tippett (1931), while the maximum is from Wilkinson (1951), who is credited with the more general approach of using an order statistic of the $p_j$.

The criterion $\min(p_1, p_2)$ leads us to reject $H_0$ if either test 1 or test 2 is strongly significant. The criterion $\max(p_1, p_2)$ is similarly seen to require at least some significance from both test 1 and test 2. Birnbaum's result opens up the possibility for many combinations between these simple types, of which Fisher's test and Stouffer's test are two prominent examples.

Graphically, we see that Fisher's combination is more sensitive to the single smallest $p$-value than Stouffer's combination is. In the Fisher test, if the first $m-1$ $p$-values already yield a test statistic exceeding the $\chi^2_{(2m)}$ significance threshold, then the $m$th test statistic cannot undo it. The Stouffer test is different. Any large but finite value of $\sum_{j=1}^{m-1} \Phi^{-1}(\widetilde{p}_j)$ can be canceled by an opposing value of $\Phi^{-1}(\widetilde{p}_m)$.

The bottom row of Figure 1 illustrates the pre-images $(\widetilde{p}_1, \widetilde{p}_2)$, where $p_j = 2\min(\widetilde{p}_j, 1-\widetilde{p}_j)$ of the regions in the top row. Those images show which combinations of one sided $\widetilde{p}$-values lead to rejection of $H_0$. Each region $\widetilde{R}$ in the bottom row of Figure 1 is symmetric with respect to replacing $\widetilde{p}_j$ by $1-\widetilde{p}_j$, and so they do not favor alternatives with concordant signs.

Figure 2 shows the rejection regions for concordant versions of the Fisher and Stouffer tests. Both of these regions devote more area to catching coordinated alternatives to $H_0$ than to split decisions. Comparing the upper left and lower right corners of these regions we see that the Stouffer version is more extreme than the Fisher test in rejecting split decisions.

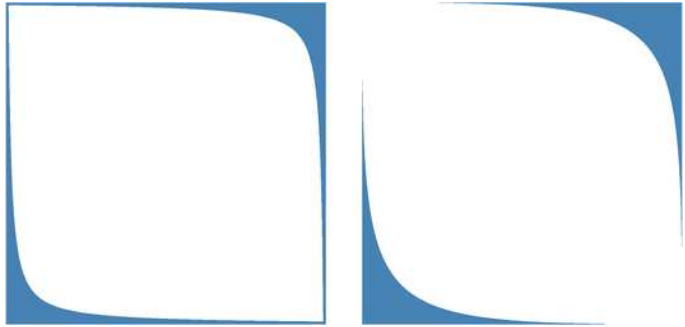

FIG. 2. *This figure shows rejection regions $\widetilde{R}$ for concordant tests, $Q^C$ and $S^C$, as described in the text. The left image shows a region for Pearson's $Q^C$ which is based on Fisher's combination. The right image shows a region for $S^C$, based on Stouffer's combination. The x and y axes in these images correspond to one sided p-values $\widetilde{p}_1$ and $\widetilde{p}_2$, rejecting $H_0$ for negative slopes at the bottom and/or left, while rejecting $H_0$ for positive slopes at the top and/or right. These tests are more sensitive to alternatives where all underlying hypothesis tests reject in the same direction than they are to split decisions. The Stouffer region extends to all corners but with a thickness that approaches zero. The Fisher region has strictly positive thickness in each corner.*

A naive reading of Condition 1 is that almost any $p$-value combination is reasonable. But some of those combinations are optimal for very unrealistic alternatives. Birnbaum (1954) goes deeper by considering alternatives in an exponential family, beginning with his second condition.

CONDITION 2. If test statistic values $(\hat{\beta}_1, \ldots, \hat{\beta}_m)$ and $(\hat{\beta}_1^*, \ldots, \hat{\beta}_m^*)$ do not lead to rejection of $H_0$, then neither does $\lambda(\hat{\beta}_1, \ldots, \hat{\beta}_m) + (1-\lambda)(\hat{\beta}_1^*, \ldots, \hat{\beta}_m^*)$ for $0 < \lambda < 1$.

Condition 2 requires that the acceptance region, in test statistic space, be convex. If the test statistics being combined are from a one parameter exponential family, then Birnbaum (1954) shows that Condition 2 is necessary for the combined test to be admissible. When the parameter space is all of $\mathbb{R}^m$, then Condition 2 is also sufficient for the combined test to be admissible. This is a consequence of the theorem in Stein (1956), Section 3. Birnbaum (1955) had this converse too, but without a condition on unboundedness of the parameter space. Matthes and Truax (1967) prove that the unboundedness is needed. Thus Condition 2 is reasonable and it rules out conjunction-based tests like the one based on $\max(p_1, \ldots, p_m)$, and more generally, all of the Wilkinson methods based on $p_{(r)}$ for $1 < r \leq m$.

Suppose that $\hat{\beta}_j \sim \mathcal{N}(\beta_j, 1)$. Then Birnbaum (1954) shows that a test which rejects $H_0$ when $\prod_{j=1}^m (1 - p_j)$ is too large, fails Condition 2. For small $\alpha$ and $m = 2$, that test has a nearly triangular rejection region $R$ including $(0,0)$ in the $(p_1, p_2)$ plane. For small $\alpha$ it would not reject at $(p_1, p_2) = (0,1)$ or even at $(0, 0.5)$. Birnbaum (1955) attributes that test to Karl Pearson through a description given by Pearson (1938). But Karl Pearson did not propose this test and it certainly is not $Q^C$. It appears that Birnbaum has misread Egon Pearson's description of Karl Pearson's test. A detailed discussion of that literature is given in Section 5.

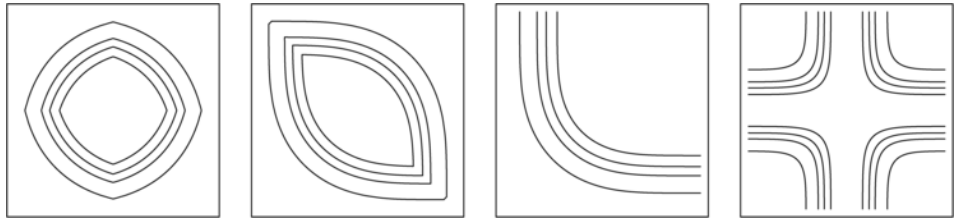

FIG. 3. *This figure shows nested decision boundaries in the space of test statistics $\hat{\beta} = (\hat{\beta}_1, \ldots, \hat{\beta}_m)$ for meta-analysis methods described in the text. We are interested in the regions where $H_0$ is not rejected. From left to right they are: Fisher's combination $Q^U$ with lozenge-shaped regions, Pearson's combination $Q^C$ with leaf-shaped regions, a left-sided combination $Q^L$ with quarter-round-shaped regions north-east of the origin and Birnbaum's version of Pearson's region, having nonconvex plus-sign-shaped regions. In each plot the significance levels are at* 0.2, 0.1, 0.05 *and* 0.01.

Theorem 1 in Section 4 shows that $Q^C$ satisfies Condition 2 for normally distributed test statistics, and so it is admissible. Figure 3 illustrates some acceptance regions for $Q^U$, $Q^C$, $Q^L$ and $\prod_{j=1}^{m}(1 - p_j)$.

In applications, Fisher's method is more widely used than Stouffer's. In a limit where sample sizes increase and test statistics become more nearly normally distributed, some methods are maximally efficient in Bahadur's sense [Bahadur (1967)]. Fisher's method is one of these, and Stouffer's is not. See Table 3 on page 44 of Hedges and Olkin (1985). Both Birnbaum (1954) and Hedges and Olkin (1985) consider Fisher's method to be a kind of default, first among equals or better.

**4. Pearson's combination.** In this section, we prove two properties of Pearson's combination $Q^C$. The acceptance regions in the second panel of Figure 3 certainly appear to be convex. We prove that his test satisfies Condition 2 (convexity), for Gaussian test statistics. Therefore, it is in fact admissible in the exponential family context for the Gaussian case. The result extends to statistics with log-concave densities. Finally, we show that the Bonferroni bound on the combination is very accurate for small combined $p$-values.

4.1. *Propagation of admissibility to $Q^C$.* We consider first the setting of Gaussian test statistics $\hat{\beta}_j \sim \mathcal{N}(\beta_j, \sigma^2/n_j)$. For simplicity, we work in terms of $t_j = \sqrt{n_j}\hat{\beta}_j/\sigma$. Under $H_0$, the $t_j$ are i.i.d. $\mathcal{N}(0,1)$. The proof that $Q^C$ is admissible uses rules that propagate convexity and quasi-convexity. It is broken into small pieces that are used for other test statistics in Sections 6 and 7.

Let $Q(t_1, \ldots, t_m)$ be a real-valued function on $\mathbb{R}^m$. The associated test rejects $H_0$ in favor of $H_A$ when $Q \geq q$. That test is admissible in the exponential family context if $\{(t_1, \ldots, t_m) | Q < q\}$ is convex and the natural parameter space is $\mathbb{R}^m$.

LEMMA 1. *For $k = 1, \ldots, n$, suppose that the test which rejects $H_0$ versus $H_A$ when $Q_k(t_1, \ldots, t_m) \geq q_k$ is admissible in the exponential family context with natural parameter space $\mathbb{R}^m$. Let the test based on $Q_k$ have level $\alpha_k$. Then the test which rejects $H_0$ when one or more of $Q_k(t_1, \ldots, t_m) \geq q_k$ hold is also admissible and has level $\alpha \leq \sum_{k=1}^n \alpha_k$.*

PROOF. Under the assumptions, the acceptance regions for the $n$ individual tests are convex. The combined test has an acceptance region equal to their intersection which is also convex. Therefore, the combined test is admissible. The level follows from the Bonferroni inequality. □

LEMMA 2. *For $k = 1, \ldots, n$, suppose that a test rejects $H_0$ versus $H_A$ when $Q_k(t_1, \ldots, t_m) \geq q_k$, where $Q_k$ is a convex function on $\mathbb{R}^m$. Then the test which rejects $H_0$ when $Q(t_1, \ldots, t_m) = \sum_{k=1}^n Q_k(t_1, \ldots, t_m) \geq q$ holds is admissible in the exponential family context with natural parameter space $\mathbb{R}^m$.*

PROOF. Under the assumptions, each $Q_k$ is a convex function, and therefore, so is their sum $Q$. Then $\{(t_1, \ldots, t_m) | Q < q\}$ is convex and so the test is admissible. □

Given a battery of admissible tests, Lemma 1 shows that we can run them all and get a new admissible test that rejects $H_0$ if any one of them rejects. Lemma 2 shows that for tests that have convex criterion functions, we can sum those functions and get a test criterion for another admissible test. Lemma 2 requires a stronger condition on the component tests $Q_k$ than Lemma 1 does. Aside from the difference between open and closed level sets, the criterion used in Lemma 1 is quasi-convexity. Sums of quasi-convex functions are not necessarily quasi-convex [Greenberg and Pierskalla (1971)].

THEOREM 1. *For $t_1, \ldots, t_m \in \mathbb{R}^m$ let*

$$Q^C = \max\left(-2\log \prod_{j=1}^m \Phi(t_j), -2\log \prod_{j=1}^m \Phi(-t_j)\right).$$

Then $\{(t_1,\ldots,t_m)|Q^C < q\}$ is convex so that Pearson's test is admissible in the exponential family context, for Gaussian data.

PROOF. The function $\Phi(t)$ is log concave. Therefore, the test that rejects $H_0$ when $Q_k = -2\log\Phi(t_k)$ is large (ignoring the other $m-1$ statistics) has a convex criterion function. Applying Lemma 2 to those tests shows that the test based on $Q^L = -2\sum_{j=1}^m \log(\Phi(t_j))$ has convex acceptance regions. The same argument starting with log concavity of $\Phi(-t)$ shows that the test based on $Q^R = -2\sum_{j=1}^m \log(\Phi(-t_j))$ has convex acceptance regions. Then Lemma 1 applied to $Q^L$ and $Q^R$ shows that $Q^C$ has convex acceptance regions. $\square$

The Gaussian case is concrete and is directly related to the motivating context. But the result in Theorem 1 holds more generally. Suppose that the density functions $(d/dt)F_j(t)$ exist and are log concave for $j = 1,\ldots,m$. Then both $F_j(t)$ and $1 - F_j(t)$ are log concave [Boyd and Vandeberghe (2004), Chapter 3.5]. Then the argument from Theorem 1 shows that the test criteria $Q^L = -2\sum_j \log(F_j(t_j))$, $Q^R = -2\sum_j \log(1-F_j(t_j))$ and $Q^C = \max(Q^L, Q^R)$ all have convex acceptance regions.

While many tests built using Lemmas 1 and 2 are admissible in the exponential family context, it may still be a challenge to find their level. The next section takes up the issue of bounding $\alpha$ for tests like $Q^C$.

4.2. *Bonferroni.* Here we show that the Bonferroni bound for $Q^C$ is very accurate, in the tails where we need it most. The Bonferroni bound for $Q^C$ is so good because it is rare for both $Q^L$ and $Q^R$ to exceed a high level under $H_0$. For notational simplicity, this section uses $X_j$ in place of $\widetilde{p}_j$.

THEOREM 2. *Let $X = (X_1,\ldots,X_m) \in (0,1)^m$ be a random vector with independent components, and put*

$$Q^L \equiv -2\log\left(\prod_{j=1}^m X_j\right), \qquad Q^R \equiv -2\log\left(\prod_{j=1}^m (1-X_j)\right)$$

*and $Q^C = \max(Q^L, Q^R)$. For $A \in \mathbb{R}$, let $P_L = \Pr(Q^L > A)$, $P_R = \Pr(Q^R > A)$ and $P_T = \Pr(Q^C > A)$. Then $P_L + P_R - P_L P_R \le P_T \le P_L + P_R$.*

COROLLARY 1. *Suppose that $X \sim U(0,1)^m$ in Theorem 2. For $A > 0$ let $\tau_A = \Pr(\chi^2_{(2m)} > A)$. Then*

$$2\tau_A - \tau_A^2 \le \Pr(Q^C > A) \le 2\tau_A.$$

PROOF. When $X \sim U(0,1)^m$, then $\Pr(Q^L > A) = \Pr(Q^R > A) = \tau_A$. The result follows by substitution in Theorem 2. □

Before proving Theorem 2, we present the concept of associated random variables, due to Esary, Proschan and Walkup (1967).

DEFINITION 2. A function $f$ on $\mathbb{R}^n$ is nondecreasing if it is nondecreasing in each of its $n$ arguments when the other $n-1$ values are held fixed.

DEFINITION 3. Let $X_1, \ldots, X_n$ be random variables with a joint distribution. These random variables are associated if $\mathrm{Cov}(f(X), g(X)) \geq 0$ holds for all nondecreasing functions $f$ and $g$ for which the covariance is defined.

LEMMA 3. *Independent random variables are associated.*

PROOF. See Section 2 of Esary, Proschan and Walkup (1967). □

LEMMA 4. *For integer $m \geq 1$, let $X = (X_1, \ldots, X_m) \in (0,1)^m$ have independent components. Set*

$$Q^L = -2\log \prod_{j=1}^m X_j \quad and \quad Q^R = -2\log \prod_{j=1}^m (1-X_j).$$

*Then for any $A_L > 0$ and $A_R > 0$,*

$$\Pr(Q^L > A_L, Q^R > A_R) \leq \Pr(Q^L > A_L)\Pr(Q^R > A_R).$$

PROOF. Let $f_1(X) = 2\log \prod_{j=1}^m X_j$ and $f_2(X) = -2\log \prod_{j=1}^m (1-X_j)$. Then both $f_1$ and $f_2$ are nondecreasing functions of $X$. The components of $X$ are independent, and hence are associated. Therefore,

$$\begin{aligned}
\Pr(Q^L &> A_L, Q^R > A_R) \\
&= \Pr(-f_1(X) > A_L, f_2(X) > A_R) \\
&= \Pr(f_1(X) < -A_L, f_2(X) > A_R) \\
&= \Pr(f_2(X) > A_R) - \Pr(f_1(X) \geq -A_L, f_2(X) > A_R) \\
&\leq \Pr(f_2(X) > A_R) - \Pr(f_1(X) \geq -A_L)\Pr(f_2(X) > A_R) \\
&= \Pr(f_2(X) > A_R)\Pr(f_1(X) < -A_L) \\
&= \Pr(Q^R > A_R)\Pr(Q^L > A_L). \qquad \square
\end{aligned}$$

In general, nondecreasing functions of associated random variables are associated. Lemma 4 is a special case of this fact, for certain indicator functions of associated variables.

PROOF OF THEOREM 2. The Bonferroni inequality yields $P_T \le P_L + P_R$. Finally, $P_T = P_L + P_R - \Pr(Q^L > A, Q^R > A)$ and $\Pr(Q^L > A, Q^R > A) \le P_L P_R$ from Lemma 4. $\square$

REMARK 1. The same proof holds for combinations of many other tests besides Fisher's. We just need the probability of simultaneous rejection to be smaller than it would be for independent tests.

**5. History of Pearson's test.** Pearson (1933) proposed the product $\prod_{i=1}^{n} F_0(X_i)$ as a way of testing whether i.i.d. observations $X_1, \ldots, X_n$ are from the distribution with CDF $F_0$. He finds the distribution of the product in terms of incomplete gamma functions and computes several examples. Pearson remarks that the test has advantages over the $\chi^2$ test of goodness of fit: small groups of observations do not have to be pooled together, and one need not approximate small binomial counts by a normal distribution. Pearson also notices that the approach is applicable more generally than testing goodness of fit for i.i.d. data, and in a note at the end, acknowledges that Fisher (1932) found the distribution earlier.

In his second paper on the topic, Pearson (1934) found a $p$-value for $\prod_j F(X_j)$ and one for $\prod_j (1 - F(X_j))$ and then advocated taking the smaller of these two as the "more stringent" test. Modern statisticians would instinctively double the smaller $p$-value, thereby applying a Bonferroni factor of 2, but Pearson did not do this.

Birnbaum [(1954), page 562] describes a test of Karl Pearson as follows:

> "Karl Pearson's method: reject $H_0$ if and only if $(1 - u_1)(1 - u_2) \cdots (1 - u_k) \ge c$, where $c$ is a predetermined constant corresponding to the desired significance level. In applications, $c$ can be computed by a direct adaptation of the method used to calculate the $c$ used in Fisher's method."

In the notation of this paper, $(1 - u_1)(1 - u_2) \cdots (1 - u_k)$ is $\prod_{j=1}^{m}(1 - p_j)$, for Figure 4 of Birnbaum (1954), and it is $\prod_{j=1}^{m}(1 - \widetilde{p}_j)$ for Figure 9. The latter (but not the former) would lead to an admissible test, had the rejection region been for small values of the product.

Birnbaum does not cite any of Karl Pearson's papers directly, but does cite Egon Pearson (1938) as a reference for Karl Pearson's test. Pearson [(1938), page 136] says,

> "Following what may be described as the intuitional line of approach, Pearson (1933) suggested as suitable test criterion one or other of the products
>
> $$Q_1 = y_1 y_2 \cdots y_n$$
>
> or
>
> $$Q'_1 = (1 - y_1)(1 - y_2) \cdots (1 - y_n).$$"

The quotation above omits two equation numbers and two footnotes but is otherwise verbatim. In the notation of this paper, $Q_1 = \prod_{j=1}^{m} \widetilde{p}_j$ and $Q_1' = \prod_{j=1}^{m}(1-\widetilde{p}_j)$. E. Pearson cites K. Pearson's 1933 paper, although it appears that he should have cited the 1934 paper instead, because the former has only $Q_1$, while the latter has $Q_1$ and $Q_1'$.

Birnbaum (1954) appears to have read E. Pearson's article as supplying two different proposals of K. Pearson, and then chose the one based on $Q_1'$, rejecting for large values of that product.

In an article published after Pearson (1933) but before Pearson (1934), David (1934), page 2, revisits the 12 numerical examples computed by Pearson (1933) and reports that in 4 of those, Pearson made a wrong guess as to which of $Q_1$ and $Q_1'$ would be smaller:

"Thus in 8 of the 12 illustrations the more stringent direction of the probability integrals was selected by mere inspection. In the other 4 cases B ought to have been taken instead of A, but in none of these four cases was the difference such as to upset the judgment as to randomness deduced from A."

Pearson (1933) computes $Q_1$ all 12 times and does not mention that this is a guess as to which product is smaller. Thus it is David's paper in which $\min(Q_1, Q_1')$ is first used (as opposed to $Q_1$). One might therefore credit David with this test, as for example, Oosterhoff (1969) does. But David credits Pearson for the method.

Birnbaum's conclusion about Pearson's test is now well established in the literature. Hedges and Olkin [(1985), page 38] write,

"Several other functions for combining $p$-values have been proposed. In 1933 Karl Pearson suggested combining $p$-values via the product

$$(1-p_1)(1-p_2)\cdots(1-p_k).$$

Other functions of the statistics $p_i^* = \text{Min}\{p_i, 1-p_i\}$, $i = 1,\ldots,k$, were suggested by David (1934) for the combination of two-sided test statistic, which treat large and small values of the $p_i$ symmetrically. Neither of these procedures has a convex acceptance region, so these procedures are not admissible for combining test statistics from the one-parameter exponential family."

The $p_i$ in this quotation might refer to either $\widetilde{p}_i$ or $p_i$ in the notation of this paper. If the former, then the product would be equivalent to $Q^R$ and would be admissible. If the latter, then the product would conform to the statistic that Birnbaum studies, but not the one Karl Pearson proposed. Furthermore, in that case the quantity $p_i^* = \min(p_i, 1-p_i)$ vanishes at $\widetilde{p}_i \in \{0, 1/2, 1\}$ and takes its maximal value at $\widetilde{p}_i \in \{1/4, 3/4\}$, which does not make sense.

**6. Power comparisons.** Admissibility of $Q^C$ is historically interesting, but for practical purposes we want to know how its power compares to other test statistics such as undirected ones, Stouffer based ones and some likelihood ratio tests developed below.

In this section, we compare the power of these tests at some alternatives of the form $(\Delta, \ldots, \Delta, 0, \ldots, 0)$ for $\Delta > 0$ where $k$ components are nonzero and $m - k$ are zero. Alternatives of the form $(-\Delta, \Delta, \ldots, \Delta, 0, \ldots, 0)$ for $\Delta > 0$ are also investigated.

Not surprisingly, concordant tests generally outperform their undirected counterparts. The most powerful method depends on the value of $k$. In applications where we know roughly how many nonzero and discordant slopes to expect, we can then identify which method will be most powerful, using the methods in this section.

The power of tests based on $S^L$, $S^R$ and $S^C$ can be handled via the normal CDF. The statistics $Q^L$, $Q^R$, $Q^U$ and $S^U$ are all sums of $m$ independent nonnegative random variables. It is therefore possible to get a good approximation to their exact distributions via the fast Fourier transform (FFT). For each of them, we use the FFT first to find their critical level (exceeded with small probability $\alpha$). The FFT is used in such a way as to get hard upper and lower bounds on the cumulative probabilities which yield hard upper and lower bounds for the critical value. Then a second FFT under the alternative hypothesis is used to compute their power. The upper limit of the power of some $Q$ comes from the upper bound on the probability of $Q$ exceeding the lower bound on its critical value. The lower limit of power is defined similarly.

For $Q^C$, the upper and lower limits on the critical value come via (2.5) applied to the bounds for $Q^L$ and $Q^R$. So do the upper and lower limits for the power.

All the computations were also done by simple Monte Carlo, with 99.9% confidence intervals. Usually the 100% intervals from the FFT were narrower than the MC intervals. But in some cases where (2.5) is not very tight, such as concordant tests at modest power, the MC intervals came out narrower.

6.1. *Alternatives to meta-analysis.* In the AGEMAP study all of the original data were present and so one could use them to form the usual test statistics instead of summing logs of $p$-values. To focus on essentials, suppose that we observe $Z_j \sim \mathcal{N}(\beta_j, 1)$ for $j = 1, \ldots, m$.

Then some very natural ways to pool the data are via $Z^R = \sum_{j=1}^m Z_j$, $Z^L = -Z^R$ and $Z^U = \sum_{j=1}^m Z_j^2$. Of course $Z^L = \sqrt{m} S^L$ and $Z^R = \sqrt{m} S^R$ but $Z^U$ is not equivalent to $S^U$. We would use $Z^R$ if the alternative were $\beta_1 = \beta_2 = \cdots = \beta_m > 0$, or even if, within the alternative $H_A$, the positive diagonal were especially important. The test criterion $Z^U$ is a more principled alternative to Fisher's $Q^U = -2 \sum_{j=1}^m \log(2\Phi(-|Z_j|))$ which does not account for the Gaussian distribution of the test statistics. Like $Q^U$ it does not favor concordant alternatives.

Marden (1985) presents the likelihood ratio test statistic (LRTS) for $H_R$ versus $H_0$. It takes the form $\Lambda^R = \sum_{j=1}^m \max(0, Z_j)^2$. The LRTS for $H_L$ versus $H_0$ is $\Lambda^L = \sum_{j=1}^{m} \max(0, -Z_j)^2$. Here we use $\Lambda^C = \max(\Lambda^L, \Lambda^R)$, a concordant version of Marden's LRTS.

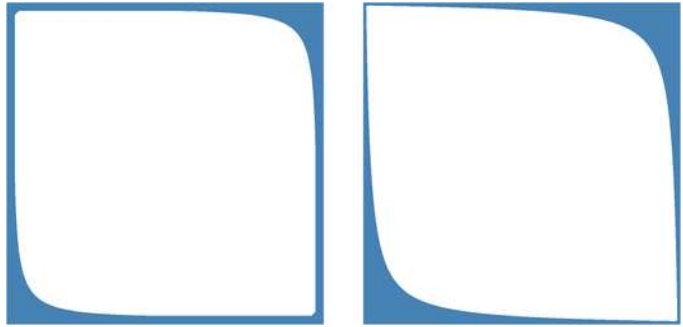

FIG. 4. *The left panel is the rejection region for $\Lambda^C$ formed by combining two Gaussian likelihood ratio tests of Marden (1985). The right panel is rejection region for Pearson's test $Q^C$. Both regions are with respect to one-tailed p-values $(\widetilde{p}_1, \widetilde{p}_2) \in (0,1)^2$.*

PROPOSITION 1. *Let $X = (X_1, \ldots, X_m) \sim U(0,1)^m$, $Z_j = \Phi^{-1}(X_j)$, and put $\Lambda^L = \sum_{j=1}^{m} \max(0, -Z_j)^2$, $\Lambda^R = \sum_{j=1}^{m} \max(0, Z_j)^2$ and $\Lambda^C = \max(\Lambda^L, \Lambda^R)$. Then $\{X | \Lambda^C \le \lambda\}$ is convex. For $A > 0$ let $\tau_A = \Pr(\chi^2_{(2m)} \ge A)$. Then $2\tau_A - \tau_A^2 \le \Pr(\Lambda^C \ge A) \le 2\tau_A$.*

PROOF. Convexity of the acceptance region for $\Lambda^C$ follows as in Theorem 1 by starting with convexity of $\max(0, \Phi^{-1}(X_j))^2$ and using Lemmas 1 and 2. For the second claim, the same argument used to prove Theorem 2 applies here, starting with nondecreasing functions $f_1(X) = -\sum_{j=1}^{m} \max(0, \Phi^{-1}(X_j))^2$ and $f_2(X) = \sum_{j=1}^{m} \max(0, -\Phi^{-1}(X_j))^2$. □

Figure 4 compares the rejection regions for $\Lambda^C$ and $Q^C$. Here we see that $Q^C$ has more of its rejection region devoted to coordinated alternatives than does $\Lambda^C$. Recalling Figure 2, we note that $S^C$ has even more of its region devoted to coordinated alternatives than $Q^C$, and so $Q^C$ is in this sense intermediate between these two tests.

Marden (1985) also presents an LRTS based on $t$-statistics. As the degrees of freedom increase the $t$-statistics rapidly approach the normal case considered here. They are not, however, an exponential family at finite degrees of freedom. If experiment $j$ gives a $t$-statistic of $T_j$ on $n_j$ degrees of freedom, then the one-sided likelihood ratio test rejects $H_0$ for large values of

$$(6.1) \quad T^R = \sum_{j=1}^{m} (n_j + 1) \log(1 + \max(T_j, 0)^2 / n_j).$$

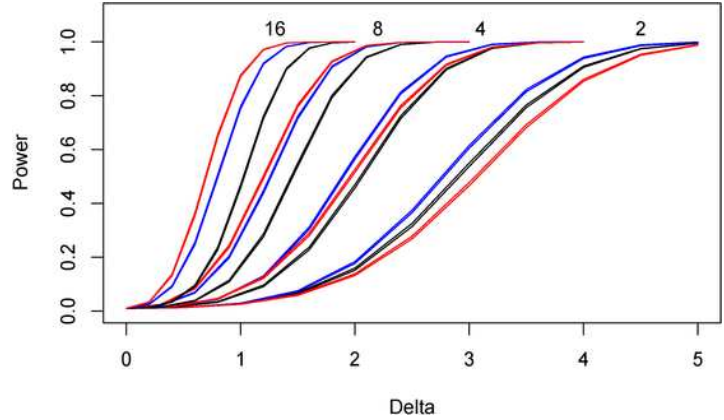

FIG. 5. *This figure shows power ranging from near* 0 *to* 1 *for a test of* $H_0$ *at level* $\alpha = 0.01$. *The alternative hypothesis is* $H_A$. *The true parameter values* $\beta$ *have* $k$ *components* $\Delta > 0$ *and* $m - k$ *components* 0. *Here* $m = 16$ *and the values of* $k$ *are printed above their power curves. The black lines are for the usual* $\chi^2$ *test statistic* $Z^U = \sum_j Z_j^2$. *The red lines are for* $Q^C$ *and the blue lines are for* $\Lambda^C$. *For each curve, upper and lower bounds are plotted, but they almost overlap.*

This is (18) of Marden (1985), after correcting a small typographical error. For large $n$ the summands in (6.1) expressed as functions of $\widetilde{p}_j \sim U(0,1)$ are very similar to $\max(0, Z_j)^2$.

6.2. *Numerical examples.* We have $H_0 : \beta = (0, \ldots, 0)$ and $H_A$ is $\beta \neq 0$. The tests will be made at level $\alpha = 0.01$. This value is a compromise between the value 0.001 used in Zahn et al. (2007) and the more conventional value 0.05. The ranking among tests is not very sensitive to $\alpha$.

Let $m = 16$ and suppose that $\beta = (\Delta, \ldots, \Delta, 0, \ldots, 0) \in H_A \subset \mathbb{R}^m$ for $\Delta > 0$. The estimates $\hat{\beta}$ are distributed as $\mathcal{N}(\beta, I_m)$. The number of nonzero components is $k \in \{2, 4, 8, 16\}$. As $\Delta$ increases the power of the tests increases. The results are shown in Figure 5.

If $k \geq 8$, then $Q^C$ (red curves) has better power than $\Lambda^C$ (blue curves), while for small $k$, $\Lambda^C$ does better. The black curves are for the test statistic $Z^U$. For $k \geq 4$, the concordant methods dominate $Z^U$. In this example, $S^C$ has the best power for $k = 16$. Power for $S^C$ is not shown but is included on later plots. Tests based on $Q^C$ do poorly in the case with $k = 2$.

Continuing this example, we now make $\Delta$ depend on $k$, so that we can get all values of $k$ onto one plot. Specifically $\Delta = \Delta_k$ is chosen so that the usual test based on $Z^U$ has power exactly 0.8. Then the best method for small $k$ arises from $\Lambda^C$, the best for the largest $k$ comes from $S^C$, while $Q^C$ is best in the middle range and is nearly best for large $k$. The central Stouffer test based on $S^U$ has power less than 0.8 over the whole range. The results are shown in Figure 6.

Finally, we consider the setting where $k - 1$ of the $\beta_j$ equal to $\Delta_k > 0$, while one of them is $-\Delta_k$. Again $\Delta_k$ is chosen so that a test based on $Z^U$ has power 0.8. Figure 7 shows the results. For small $k$ $\Lambda^C$ works best, while for larger $k$, $Q^C$ works best. The Stouffer test $S^C$ is best when $k = 16$, but loses power quickly as $k$ decreases.

An online supplement at http://stat.stanford.edu/~owen/reports/KPearsonSupplement contains 32 plots like the figures shown here. These

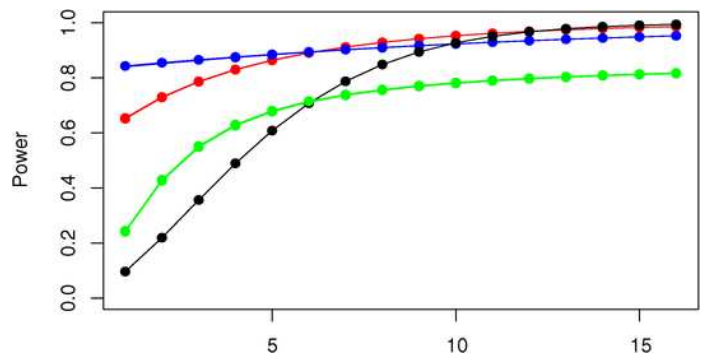

FIG. 6. *This figure shows the power of various tests of* $H_0 : \beta = 0$ *when* $\beta = (\Delta, \ldots, \Delta, 0, \ldots, 0) \in \mathbb{R}^m$. *The number $k$ of nonzero components ranges from 1 to $m = 16$ on the horizontal axis. For each $k$, $\Delta$ was chosen so that a test based on* $Z^U = \sum_{j=1}^m Z_j^2$ *has power 0.8 of rejecting $H_0$. Results are given for $Q^C$ (red), $\Lambda^C$ (blue), $S^C$ (black) and $S^U$ (green).*

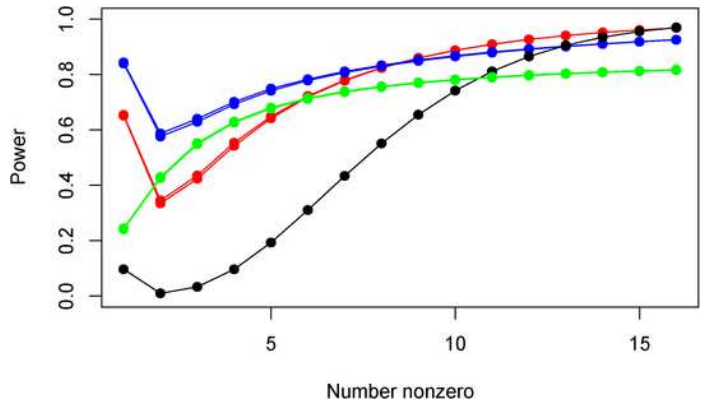

Number nonzero

FIG. 7. *This figure shows the power of various tests of* $H_0 : \beta = 0$ *when* $\beta = (-\Delta, \Delta, \ldots, \Delta, 0, \ldots, 0) \in \mathbb{R}^m$. *The number $k$ of nonzero components ranges from 1 to $m = 16$ on the horizontal axis. There is always one negative component. For each $k$, $\Delta$ was chosen so that a test based on* $Z^U = \sum_{j=1}^m Z_j^2$ *has power 0.8 of rejecting $H_0$. Results are given for $Q^C$ (red), $\Lambda^C$ (blue), $S^C$ (black) and $S^U$ (green). In this figure Monte Carlo was used for $Q^C$ and $\Lambda^C$.*

figures are also given in the technical report Owen (2009). The cases considered have $m \in \{4, 8, 12, 16\}$, $\alpha \in \{0.01, 0.05\}$ and power in $\{0.8, 0.5\}$. The number $k$ of nonzero components ranges from 1 to $m$. In one set of plots there are $k$ positive components, while the other has $k-1$ positive and 1 negative component.

Those other figures show the same patterns as the one highlighted here. Among the concordant tests, $\Lambda^C$ is best for small $k$, then $Q^C$ is best for medium sized $k$ and finally $S^C$ is best for the largest $k$. When there is one negative component, $S^C$ is most adversely affected, and $\Lambda^C$ least. That is, the tests that gain the most from coordinated alternatives, lose the most from a discordant component. The case $k=2$ is hardest when there is one negative component, for then $\beta$ contains two nonzero components of opposite sign. The Stouffer test $S^U$ is never competitive with $Z^U$, though the gap is small when $k = m$.

When $2 \leq k \leq m = 4$ and one component is negative, then $\beta$ does not really have concordant signs, and in these cases $Q^C$, $\Lambda^C$ and $S^C$ have less power than $S^U$ which has less power than $Z^U$.

6.3. *The original application.* The AGEMAP study of Zahn et al. (2007) used $Q^C$ to filter genes. That setting featured 8932 genes to be tested in $m = 16$ tissues with data from $n = 40$ animals (slightly fewer when some data were missing). There were 5 male and 5 female animals at each of 4 ages. For gene $i$ and tissue $j$ there was a regression of 40 gene expression values on age, sex and an intercept. The $p$-value for gene $i$ and tissue $j$ was based on a $t$-test for the age coefficient, usually with 37 degrees of freedom.

Of the 16 tissues investigated, only 9 appeared to show any aging, and so much of the analysis was done for just those 9.

The biological context made it impossible to specify a single alternative hypothesis, such as $H_L$, $H_R$ or $H_L \cup H_R$ to use for all genes. Instead it was necessary to screen for interesting genes without complete knowledge of which tissues would behave similarly. Also, the fates of tissues differ in ways that were not all known beforehand. One relevant difference is that the thymus involutes (becomes fatty) while some other tissues become more fibrous. There are likely to be other differences yet to be discovered that could make one tissue age differently from others.

The 8932 genes on any given microarray in the study were correlated. There were also (small) correlations between genes measured in two different tissues for the same animal. It is, of course, awkward to model the correlations among an $8932 \times 16$ matrix of observations based on a sample of only 40 such matrices. An extensive simulation was conducted in Owen (2007). In that simulation, genes were given known slopes and then errors were constructed by resampling residual matrices for the 40 animals. By resampling residuals, some correlations among genes and tissues were retained. In each Monte Carlo sample, the genes were ranked by all the tests, and ROC curves described which test method was most accurate at ranking the truly age-related genes ahead of the others.

In that Monte Carlo simulation, there were $k \in \{3, 5, 7\}$ nonzero slopes out of $m \in \{9, 16\}$. Various values of $\Delta$ were used. The Fisher-based tests consistently outperformed Stouffer-style tests, as statistical theory would predict. The concordant tests usually outperformed central tests, and were almost as good as the one-sided tests that we would use if we knew the common sign of the nonzero $\beta_j$. The exceptions were for cases with $k = 3$ and $m = 16$ and large $\Delta$. Then the $Q^U$ tests could beat the $Q^C$ tests. For $k = 3$ and $m = 16$ and small $\Delta$, the undirected and concordant tests were close. The left-hand side of Figure 6 confirms that pattern.

6.4. *Power conclusions.* The methods of this section can be used to compare different combinations of tests. Given precise information, such as a prior distribution $\pi(\beta)$ for nonnull $\beta$, one could employ a weighted sum of power calculations to estimate $\int_{\mathbb{R}^m} \pi(\beta) \Pr(Q \ge Q^{1-\alpha} | \beta) \, d\beta$ for each test $Q$ in a set of candidates.

Given less precise qualitative information, that the alternatives are likely to be concordant, we can still make a rough ordering of the methods when we have some idea how many concordant nonnull hypotheses there may be. Of the concordant methods compared above, the LRT $\Lambda^C$ is best where there are a few concordant nonnulls, then Pearson's $Q^C$ is best if we expect mostly concordant nonnulls, and finally Stouffer's method $S^C$ came out best only when the concordant nonnulls were unanimous or nearly so.

**7. Recent literature.** Although combination of $p$-values is a very old subject, there seems to be a revival of interest. Here, a few works related to the present setup are outlined.

Whitlock (2005) takes the strong view that any discordant test means that the null hypothesis should not be rejected. He gives the example of inbred animals being significantly larger than outbred in one study but significantly smaller in another study, and says that the results should then cancel each other out. In some contexts, such as AGEMAP, such cancellation would not make sense. Whitlock (2005) has a strong preference for the Stouffer test over Fisher style of tests. He reports better performance for Stouffer-style tests in a simulation, which is not what statistical theory would predict. His computations, however, are not equivalent to the Fisher statistics reported here. For example, he reports that,

"A result of $P = 0.99$ is as suggestive of a true effect as is a result of $P = 0.01$. Yet with Fisher's test, if we get two studies with $P = 0.01$, the combined $P$ is 0.001, but with two studies with $P = 0.99$ the combined $P$ is 0.9998. One minus 0.9998 is 0.0002. The high $P$ and low $P$ results differ by an order of magnitude, yet the answer should be the same in both cases."

Interpreting the $P$-values as one-tailed $\widetilde{p}_j$-values, Fisher's test $Q^U$ uses $p_j = 0.02$ in both cases for a combined $p$-value of $\Pr(\chi^2_{(4)} \geq -2\log(0.02^2)) \doteq 0.0035$. It is reasonable to conclude that the poor performance of Fisher's test reported in Whitlock (2005) does not apply to $Q^U$.

Benjamini and Heller (2007) consider "partial conjunction" alternative hypotheses,

(7.1)  $H_r : \beta_j \neq 0$  for at least $r$ values of $j \in \{1, \ldots, m\}$.

In that setting, one decides a priori that there must be at least $r$ false nulls among the $m$ component hypotheses for the finding to be useful.

They present several tests based on combinations of all but the $r - 1$ most significant of the tests. Their combination methods include Fisher's and Stouffer's as well as one due to Simes (1986). The partial conjunction null described above can be applied using $\widetilde{p}_j$, $1 - \widetilde{p}_j$, or $p_j$ to get left, right and undirected versions. A concordant version could also be useful if it were only of interest to find hypotheses rejected in the same direction by at least $r$ tests. To get a concordant version, one takes twice the smaller of the left- and right-sided combined $p$-values.

For $r > 1$, nontrivial acceptance regions based on combinations of only the *least significant* $m - r + 1$ $p$-values are not convex because they include any point on any of the $m$ coordinate axes in $\mathbb{R}^d$. As a result, the tests are not admissible versus the point hypothesis $H_0$ in the exponential family context. The alternatives $H_r$ in (7.1) for $r > 1$ are not simple null hypotheses though, and so the tests may be admissible for their intended use.

**8. Discussion.** This paper has shown that Pearson's method $Q^C$ really is admissible, and perhaps surprisingly, is competitive with standard methods based on the raw data, not just the $p$-values. The context where it is competitive is one where the truly nonzero components of the parameter vector are predominantly of one sign. We have also studied a concordant LRT test $\Lambda^C$ which performs well when the number of concordant alternatives is slightly less.

Also a very simple Bonferroni calculation proved to be very accurate for finding critical values of tests. It is less accurate for computing modest power levels.

In a screening setting like AGEMAP, we anticipate that noise artifacts could give rise to values $\hat{\beta}_j$ with arbitrary patterns of signs, while the true biology is likely to be dominated by concordant signs. In early stages of the investigation, false discoveries are considered more costly than false nondiscoveries because the former lead to lost effort. Later when the aging process is better understood, there may be greater value in finding those genes that are strongly discordant. In that case, combination statistics which favor discordant alternatives may be preferred.

Finally, this work has uncovered numerous errors in earlier papers. I do not mean to leave the impression that the earlier workers were not careful, either in an absolute or relative sense. The subject matter is very slippery.

## APPENDIX: COMPUTATION

We want to get the $1-\alpha$ quantile of the distribution of $Q = \sum_{j=1}^{m} Y_j$ where $Y_j$ are independent but not necessarily identically distributed random variables on $[0, \infty)$. The case of random variables with a different lower bound, possibly $-\infty$, is considered in a remark below. We suppose that $Y_j$ has cumulative distribution function $F_j(y) = \Pr(Y_j \le y)$ which we can compute for any value $y \in [0, \infty)$.

**A.1. Convolutions and stochastic bounds.** Because the $Y_j$ are independent, we may use convolutions to get the distribution of $Q$. Convolutions may be computed rapidly using the fast Fourier transform. A very fast and scalable FFT is described in Frigo and Johnson (2005) who make their source code available. Their FFT on $N$ points is tuned for highly composite values of $N$ (not just powers of 2) while costing at most $O(N \log(N))$ time even for prime numbers $N$. Thus one does not need to pad the input sequences with zeros.

There are several ways to apply convolutions to this problem. For a discussion of statistical applications of the FFT, including convolutions of distributions, see the monograph by Monahan (2001). The best-known method convolves the characteristic functions of the $Y_j$ to get that of $Q$ and then inverts that convolution. But that method brings aliasing problems. We prefer to convolve probability mass functions. This replaces the aliasing problem by an overflow problem that is easier to account for.

We write $F \otimes G$ for the convolution of distribution functions $F$ and $G$. Our goal is to approximate $F_Q = \bigotimes_{j=1}^{m} F_j$. We do this by bounding each $F_j$ between a stochastically smaller discrete CDF $F_j^-$ and a stochastically larger one $F_j^+$, both defined below. Write $F_j^- \preccurlyeq F_j \preccurlyeq F_j^+$ for these stochastic inequalities. Then from

$$(A.1) \qquad \bigotimes_{j=1}^{m} F_j^- \preccurlyeq F_Q \preccurlyeq \bigotimes_{j=1}^{m} F_j^+,$$

we can derive upper and lower limits for $\Pr(Q \ge Q^*)$ for any specific value of $Q^*$.

The support sets of $F_j^-$ and $F_j^+$ are

$$S_{\eta,N} = \{0, \eta, 2\eta, \ldots, (N-1)\eta\} \quad \text{and} \quad S_{\eta,N}^+ = S_{\eta,N} \cup \{\infty\},$$

respectively, for $\eta > 0$. The upper limit has

$$F_j^+(y) = \begin{cases} F_j(y), & y \in S_{\eta,N}, \\ 1, & y = \infty. \end{cases}$$

In the upper limit, any mass between the points of $S_{\eta,N}^+$ is pushed to the right. For the lower limit, we push mass to the left. If $F_j$ has no atoms in $S_{\eta,N}$, then

$$F_j^-(y) = \begin{cases} F_j(y + \eta), & y/\eta \in \{0, 1, 2, \ldots, N-2\}, \\ 1, & y = (N-1)\eta, \end{cases}$$

and otherwise we use $\lim_{z \downarrow (y+\eta)} F_j(z)$ for the first $N-1$ levels. We do not need to put mass at $-\infty$ in $F_j^-$ because $F_j$ has support on $[0, \infty)$.

It should be more accurate to represent each $F_j$ at values $(i + 1/2)\eta$ for $0 \leq i < N$ and convolve those approximations [see Monahan (2001)]. But that approach does not give hard upper and lower bounds for $F_Q$.

Suppose that $F$ and $G$ both have support $S_{\eta,N}^+$ with probability mass functions $f$ and $g$ respectively. Then their convolution has support $S_{\eta,2N-1}^+$. The mass at $\infty$ in the convolution is $(f \otimes g)(\infty) = f(\infty) + g(\infty) - f(\infty)g(\infty)$. The mass at multiples 0 through $2N - 2$ times $\eta$, is the ordinary convolution of mass functions $f$ and $g$,

$$(f \otimes g)(k\eta) = \sum_{i=0}^{k} f(i\eta)g((k-i)\eta).$$

The CDF $F \otimes G$ can then be obtained from the mass function $f \otimes g$. Thus the convolutions in (A.1) can all be computed by FFT with some additional bookkeeping to account for the atom at $+\infty$.

When $F$ and $G$ have probability stored at $N$ consecutive integer multiples of $\eta > 0$, then their convolution requires $2N - 1$ such values. As a result, the bounds in (A.1) require almost $mN$ storage. If we have storage for only $N$ finite atoms the convolution could overflow it. We can save storage by truncating the CDF to support $S_{\eta,N}^+$ taking care to round up when convolving factors of the upper bound and to round down when convolving factors of the lower bound.

For a CDF $F$ with support $S_{\eta,M}^+$ where $M \geq N$, define $\lceil F \rceil_N$ with support $S_{\eta,N}^+$ by

$$\lceil F \rceil_N(i\eta) = F(i\eta), \quad 0 \leq i < N, \quad \text{and} \quad \lceil F \rceil_N(\infty) = 1.$$

That is, when rounding $F$ up to $\lceil F \rceil_N$, all the atoms of probability on $N\eta$ through $(M-1)\eta$ inclusive are added to the atom at $+\infty$.

To round this $F$ down to $S_{\eta,N}$, we may take

$$\lfloor F \rfloor_N(i\eta) = F(i\eta), \quad 0 \leq i < N-1 \quad \text{and} \quad \lfloor F \rfloor_N((N-1)\eta) = 1.$$

When rounding $F$ down to $\lfloor F \rfloor_N$, all the atoms on $N\eta$ through $(M-1)\eta$ and $+\infty$ are added to the atom at $F((N-1)\eta)$. This form of rounding never leaves an atom at $+\infty$ in the stochastic lower bound for a CDF. It is appropriate if the CDF being bounded is known to be proper. If the CDF to be bounded might possibly be improper with an atom at $+\infty$, then we could instead move only the atoms of $F$ on $N\eta$ through $(M-1)\eta$ to $\lfloor F \rfloor_N((N-1)\eta)$, leave some mass at $\infty$, and get a more accurate lower bound.

The upper and lower bounds for $F_Q$ are now $F_Q^{m+}$ and $F_Q^{m-}$, where

$$F_Q^{j+} = \lceil F_Q^{(j-1)+} \otimes F_j^+ \rceil_N, \qquad j=1,\ldots,m,$$
$$F_Q^{j-} = \lfloor F_Q^{(j-1)-} \otimes F_j^- \rfloor_N, \qquad j=1,\ldots,m,$$

and $F_Q^{0\pm}$ is the CDF of a point mass at 0.

If all of the $F_j$ are the same then one may speed things up further by computing $F_1^{2^r+}$ via $r-1$ FFTs in a repeated squaring sequence, and similarly for $F_1^{2^r-}$. For large $m$, only $O(\log(m))$ FFTs need be done to compute $F_1^{m\pm}$.

REMARK 2. If one of the $F_j$ has some support in $(-\infty, 0)$ then changes are required. If $F_j$ has support $[-A_j, \infty)$ for some $A_j < \infty$ then we can work with the random variable $Y_j + A_j$ which has support $[0, \infty)$. The convolution of $F_j$ and $F_k$ then has support starting at $-(A_j + A_k)$. If $F_j$ does not have a hard lower limit like $A_j$ then we may adjoin an atom at $-\infty$ to the CDF representing its stochastic lower bound. As long as we never convolve a distribution with an atom at $+\infty$ with another distribution having an atom at $-\infty$, the results are well defined CDFs of extended real-valued random variables.

**A.2. Alternative hypotheses.** In this section, we get expressions for the CDFs $F_j$ that need to be convolved. We suppose that $\hat{\beta}_j$ are independent $\mathcal{N}(\beta_j, 1)$ random variables for $j=1,\ldots,m$. The null hypothesis is $H_0\colon \beta = 0$ for $\beta = (\beta_1, \ldots, \beta_m)$.

The left, right and undirected test statistics take the form $\sum_{j=1}^m Y_j$, where $Y_j = t(\hat{\beta}_j)$ for a function $t$ mapping $\mathbb{R}$ to $[0, \infty)$. Large values of $Y_j$ represent stronger evidence against $H_0$. The concordant test statistics are based on the larger of the left- and right-sided sums.

The likelihood ratio tests $\Lambda^L$, $\Lambda^R$ and $\Lambda^U$ are sums of

$$Y_{Lj} = \max(-\hat{\beta}_j, 0)^2, \qquad Y_{Rj} = \max(\hat{\beta}_j, 0)^2 \quad \text{and} \quad Y_{Uj} = \hat{\beta}_j^2,$$

respectively. After elementary manipulations, we find $F_j$ for these tests via

$$\Pr(Y_{Lj} \leq y) = \Phi(\sqrt{y} + \beta_j), \qquad \Pr(Y_{Rj} \leq y) = \Phi(\sqrt{y} - \beta_j)$$

and
$$\Pr(Y_{Uj} \leq y) = \Phi(\sqrt{y} - \beta_j) - \Phi(-\sqrt{y} - \beta_j).$$

The Fisher test statistics $Q^L$, $Q^R$ and $Q^U$ are sums of
$$Y_{Lj} = -2\log(\Phi(\hat{\beta}_j)), \qquad Y_{Rj} = -2\log(\Phi(-\hat{\beta}_j))$$

and
$$Y_{Uj} = -2\log(2\Phi(-|\hat{\beta}_j|)),$$

respectively. The corresponding $F_j$ are given by
$$\Pr(Y_{Lj} \leq y) = \Phi(\beta_j - \Phi^{-1}(e^{-y/2})),$$
$$\Pr(Y_{Rj} \leq y) = \Phi(-\beta_j - \Phi^{-1}(e^{-y/2}))$$

and
$$\Pr(Y_{Uj} \leq y) = \Phi(\beta_j - \Phi^{-1}(\tfrac{1}{2}e^{-y/2})) - \Phi(\beta_j + \Phi^{-1}(\tfrac{1}{2}e^{-y/2})).$$

For three of the Stouffer statistics, no FFT is required because $S^R \sim \mathcal{N}(m^{-1/2} \times \sum_{j=1}^{m} \beta_j, 1)$, $S^L = -S^R$ and $S^C = |S^R|$. The remaining Stouffer statistic $S^U$ is the sum of $Y_{Uj} = |\hat{\beta}_j|/\sqrt{m}$, with
$$\Pr(Y_{Uj} \leq y) = \Phi(\sqrt{m}y - \beta_j) - \Phi(-\sqrt{m}y - \beta_j).$$

**A.3. Convolutions for power calculations.** The computations for this paper were done via convolution using $N = 100{,}000$ and $\eta = 0.001$. Some advantage might be gained by tuning $N$ and $\eta$ to each case, but this was not necessary. The convolution approach allows hard upper and lower bounds for probabilities of the form $\Pr(Q \leq Q^*)$ for given distributions $F_j$. For practical values of $N$, the width of these bounds is dominated by discretization error in approximating $F$ at $N$ points. Empirically, it decreases like $O(1/N)$, for a given $m$, as we would expect because apart from some mass going to $+\infty$, any mass being swept left or right moves at most $O(\eta/N)$. For very large values of $N$, the numerical error in approximating $\Phi^{-1}$ would become a factor.

Each convolution was coupled with a Monte Carlo computation of $N$ sample realizations. Partly this was done to provide a check on the convolutions, but in some instances, the Monte Carlo answers were more accurate.

The possibility for Monte Carlo to be sharper than convolution arises for test statistics like $Q^C = \max(Q^R, Q^L)$. Suppose that $Q^L$ and $Q^R$ are negatively associated and that we have bounds $F_Q^{L-} \preccurlyeq F_Q^L \preccurlyeq F_Q^{L+}$ and $F_Q^{R-} \preccurlyeq F_Q^R \preccurlyeq F_Q^{R+}$. Even if $F_Q^{L-} = F_Q^{L+}$ and $F_Q^{R-} = F_Q^{R+}$, we still do not get arbitrarily narrow bounds for $F_Q^C$. In particular, increasing $N$ will not suffice to get an indefinitely small interval.

When $Q^L$ and $Q^R$ are negatively associated, we can then derive from Theorem 2 that $F_Q^{C-} \preccurlyeq F_Q^C \preccurlyeq F_Q^{C+}$ where for $i = 0, \ldots, N-1$,

$$F_Q^{C-}(i\eta) = F_Q^{R-}(i\eta)F_Q^{L-}(i\eta),$$
$$F_Q^{C+}(i\eta) = \max(0, F_Q^{R+}(i\eta) + F_Q^{L+}(i\eta) - 1)$$

and $F_Q^{C+}(\infty) = 1$.

Surprisingly, this is often enough to get a very sharp bound on $Q^C$. But in some other cases, the Monte Carlo bounds are sharper.

Monte Carlo confidence intervals for $\Pr(Q^C > Q^*)$ were computed by a formula for binomial confidence intervals in Agresti (2002), page 15. This formula is the usual Wald interval after adding pseudo-counts to the number of successes and failures. For a $100(1-\alpha)\%$ confidence interval one uses

$$\widetilde{\pi} \pm \Phi^{-1}(1-\alpha/2)\sqrt{\widetilde{\pi}(1-\widetilde{\pi})/N},$$

where $\widetilde{\pi} = (N\hat{\pi} + \Phi^{-1}(1-\alpha/2)^2/2)/(N + \Phi^{-1}(1-\alpha/2)^2)$, and $\hat{\pi}$ is simply the fraction of times that $Q^C > Q^*$ was observed in $N$ trials. For $\alpha = 0.001$, this amounts to adding $\Phi^{-1}(0.9995)^2 \doteq 10.8$ pseudo-counts split equally between successes and failures, to the $N$ observed counts.

**Acknowledgments.** I thank Stuart Kim and Jacob Zahn of Stanford's department of Developmental Biology for many discussions about $p$-values and microarrays that shaped the ideas presented here. I thank Ingram Olkin for discussions about meta-analysis and John Marden for comments on the likelihood ratio tests for $t$-distributed data.

I thank two anonymous reviewers for their comments, and especially for mentioning the article by Stein.

Department of Statistics
Sequoia Hall
390 Serra Mall
Stanford, California 94305
USA
E-mail: owen@stat.stanford.edu
URL: http://stat.stanford.edu/~owen